\pdfoutput=1
\documentclass[10pt]{article}%,letterpaper

\usepackage{selinput}
\makeatletter
\usepackage{amscd}
\usepackage{amsmath}
\usepackage{latexsym}
\usepackage{amsfonts}
\usepackage{amssymb}
\usepackage{amsthm}
\usepackage{cite}
\usepackage[%pdfstartview=FitH,
            CJKbookmarks=true,
            bookmarksnumbered=true,
            bookmarksopen=true,
            colorlinks=true, %注释掉此项则交叉引用为彩色边框(将colorlinks和pdfborder 同时注释掉)
            citecolor=magenta,% magenta , cyan
            linkcolor=blue,
            ]{hyperref}       % hyperref 宏包通常要求放在导言区的最后!!!

\newtheorem{thm}{\textbf Theorem}[section]
\newtheorem{lem}{\textbf Lemma}[section]
\newtheorem{rem}{\textbf Remark}[section]

\newcommand{\md}{\mbox{d}}
\newcommand{\mr}{\mathbb{R}}

\newcommand{\pa}{\partial}

\begin{document}

\begin{titlepage}
\title{\bf  A regularity criterion  for the 3D Boussinesq equations
 in homogeneous
Besov spaces with negative  indices  }
\author{ Mianlu Zou$^{1}$, Qiang Li $^{2}$\thanks{Corresponding Author: Q. Li}
        \\ $^{1}$ School of Big Data and Artificial Intelligence,
        Xinyang College, Xinyang, 464000, China\\
        \\ $^{2}$ School of Mathematics and Statistics,
         Xinyang College, Xinyang, 464000, China\\
       ( zoumianlu@163.com, liqiang5412@163.com)
          }
\date{}
\end{titlepage}
\maketitle

\begin{abstract}
    In this paper, we study the regularity criteria for the 3D Boussinesq equations in
terms of one partial derivative of the velocity in Besov spaces. More precisely, it is proved that  if the   velocity  $u$ holds
$$\int_{0}^{T}\| \partial_{3} u\|_{\dot{B}_{\infty,\infty}^{-r}}^{\frac{2}{1-r}}\md t<\infty,\ with\  \ 0\leq r<1,$$
   then the  solution $(u, \theta)$ is regular on $[0,T]$.
\end{abstract}

\vspace{.2in} {\bf Key words:}\quad Boussinesq equations; regularity criterion; Besov spaces

\vspace{.2in} {\bf Mathematics Subject Classification:}\quad 35Q35; 35B65; 76D03

\vspace{.2in}

%%\newpage

%%%%%%%%%%%%%%%%%%%%%%%%%%%%%%%%%%%%%%%%%%%%%%%%%%%%%%%%%%%%%%%%%%%%%%%%%%%%%%%%%%%%%%%%%%%%%%%%%%%%%%%%%%%%%%%%%%%%%%%%%%%%

\section{Introduction}
\setcounter{equation}{0}
\vskip .1in
The Boussinesq equations (\ref{1.1}) is  widely used to model the oceanic and atmospheric motions and study atmospheric sciences (see e.g.\cite{A-2003}).
In this paper, we  investigate  the regularity problem for the following three-dimensional incompressible Boussinesq equations:
\begin{align}\label{1.1}
 \begin{cases}
\partial_t u + u\cdot \nabla u - \mu\Delta u+\nabla p = \theta e_{3},\\
\partial_t \theta + u\cdot \nabla \theta - \nu \Delta\theta=0, \\
\nabla \cdot u =0,\\
u(x, 0) = u_0(x), \theta(x, 0) = \theta_0(x),
\end{cases}
 \end{align}
where $ u=(u_{1}, u_{2}, u_{3})\in\mathbb{R}^{3} $ is  the velocity vector field, $ \theta\in\mathbb{R} $ is the scalar function temperature, and $ p\in\mathbb{R} $ represents the scalar function pressure.  And $e_{3}=(0,0,1)$, the parameters $\mu$ and $\nu$ are positive constants, for simplicity, we shall assume that $\mu=\nu=1 $ in this paper.

It should be pointed out that the system (\ref{1.1})  reduces to the 3D classical incompressible
Navier-Stokes equations  when $\theta=0$. Therefore, like  Navier-Stokes equations, it is necessary and meaningful to study the regularity problem of the weak solutions for system (\ref{1.1}). Up to now, many  regularity criteria for system (\ref{1.1}) have been proposed, for references see \cite{Y-2012,Fan-Z,N-1999,S-2020, H-2010,H-2011, F-2020,F-2021,Ye-2016} and references therein.
It is interesting to investigate the regularity criteria that are formed of derivative in  vertical direction on velocity, which firstly arose for Navier-Stokes equations. In \cite{Cao-2010}, Cao  proved that if the following condition holds
 \begin{align}\label{1.2}
\int_{0}^{T}\| \pa_{3}u\|_{L^{\alpha}}^{\beta}\md t<\infty,\ with \ \frac{3}{\alpha}+\frac{2}{\beta}\leq2\ , \alpha\geq 1,
\end{align}
then the solution $u$ becomes regular. Later, this result had been extended to other fluid equations (see\cite{C-2010,J-Z}).  Recently, Wu  \cite{F-2020}  showed a regularity criterion for Boussinesq equations. Namely,  if the velocity $u$ satisfies
\begin{align}\label{1.3}
\int_{0}^{T}\| \pa_{3}u\|_{\dot{B}_{\infty,\infty}^{0}}^{2}\md t<\infty,
\end{align}
 then the solution $u$ is regular on $(0.T]$. Here $\dot{B}_{\infty,\infty}^{0}$ is the homogeneous Besov space, for the definition and property refer to \cite{H-2011}.

 The purpose of this paper is to extend the   condition (\ref{1.3}) to the homogeneous Besov space with negative  indices. our main result is stated as follows:

\begin{thm}\label{thm1}
 Assume the initial data $(u_{0}, \theta_{0})\in H^{1}(\mathbb{R}^{3})\times H^{1}(\mathbb{R}^{3})$, and  $(u,\theta)$ be a smooth solution of (\ref{1.1}) on $(0,T)$ for  $T>0$. If velocity $u$  satisfis
\begin{align}\label{1.4}
\int_{0}^{T}\| \partial_{3} u\|_{\dot{B}_{\infty, \infty}^{-r}}^{\frac{2}{1-r}}\md t<\infty,\ with\  \ 0\leq r<1,
\end{align}
then the  solution $(u,\theta)$ is regular on $(0,T]$.
\end{thm}
\begin{rem}Noting that  (\ref{1.4}) reduces to  (\ref{1.3}) when $r=0$, thus the regularity criterion (\ref{1.4}) contains the regularity criterion (\ref{1.3}).
\end{rem}

\vskip .2in

\section{Proof of Theorem 1.1}

\setcounter{equation}{0}
This section is devoted to prove Theorem \ref{thm1}. Before the proof, we give
the following  lemmas  which will play
an important role in our discussion.
\begin{lem}\label{lemma 2.1}

 (Page 82 in \cite{H-J-2011}). Let $1<q<p<\infty$ and $\alpha$ be a positive real number. Then there exists a constant
C such that
\begin{align}\label{2.1}
\|f\|_{L^{p}}\leq C\|f\|_{\dot{B}_{\infty,\infty}^{-\alpha}}^{1-\theta}\| f\|_{\dot{B}_{q,q}^{\beta}}^{\theta}, with\ \beta=\alpha(\frac{p}{q}-1), \theta=\frac{q}{p}.
\end{align}
\end{lem}

\begin{lem}\label{lemma 2.2}
 (Page 2 in \cite{C-2010}) Let  $f\in H^{1}(\mr^{3})$, then exists a constant C such that
\begin{align}\label{2.2}
 \|f\|_{L^{p}}\leq C\|\pa_{1}f\|_{L^{p_{1}}}^{\frac{1}{3}}\|\pa_{2}f\|_{L^{p_{2}}}^{\frac{1}{3}}
 \|\pa_{3} f\|_{L^{p_{3}}}^{\frac{1}{3}},
 \end{align}
 where $1\leq p_{1},p_{2},p_{3}<\infty, 1+\frac{3}{p}=\frac{1}{p_{1}}+\frac{1}{p_{2}}+\frac{1}{p_{3}}$.
 \end{lem}

 \begin{lem}\label{lemma 2.3}

 Assume the initial data $(u_{0}, \theta_{0})\in H^{1}(\mathbb{R}^{3})\times H^{1}(\mathbb{R}^{3})$, $u$ and $\theta$ is a pair smooth solution of the system (\ref{1.1}). If  the following condition is satisfied
\begin{align}\label{2.3}
\int_{0}^{T}\| \pa_{3} u\|_{\dot{B}_{\infty, \infty}^{-r}}^{\frac{2}{1-r}}\md t<\infty,\ with\  \ 0\leq r<1,
\end{align}
then we have
   $$\|\pa_{3} u\|_{L^{2}}^{2}+\int_{0}^{T}\|\pa_{3}\nabla u\|_{L^{2}}^{2}\md t<\infty.$$

  \end{lem}
\textbf{Proof:}
Firstly, it is not difficult to get the $L^{2}$ estimate of $u$ and $\theta$.  Taking the $L^{2}$ inner product with $u$ and $\theta$ to the equations $(1.1)_{1}$ and $(1.1)_{2}$ respectively, one gets
\begin{align*}
&\frac{1}{2}\frac{\md}{\md t}(\|u\|_{L^{2}}^{2}+\|\theta\|_{L^{2}}^{2})+\|\nabla u\|_{L^{2}}^{2}+\|\nabla \theta\|_{L^{2}}^{2}
=-\int_{\mr^{3}}\theta e_{3}u\md x
\\
&\leq\|\theta\|_{L^{2}}\|u\|_{L^{2}}
\leq C(\|u\|_{L^{2}}^{2}+\|\theta\|_{L^{2}}^{2}),
\end{align*}
thus
\begin{align*}
(\|u\|_{L^{2}}^{2}+\|\theta\|_{L^{2}}^{2})+\int_{0}^{T}(\|\nabla u\|_{L^{2}}^{2}+\|\nabla \theta\|_{L^{2}}^{2})\md t
\leq C.
\end{align*}
Applying $\partial_{3}$ to the first equation in $(1.1)$, and then  taking the $L^{2}$ inner product  with $\partial_{3}u$ for the resulting equation,  we have
 \begin{align}\label{2.4}
\frac{1}{2}\frac{\md}{\md t}\|\pa_{3} u\|_{L^{2}}^{2}+\|\pa_{3}\nabla u\|_{L^{2}}^{2}&=
-\int_{\mr^{3}}\pa_{3} (u\cdot\nabla u)\cdot\pa_{3} u\md x+\int_{\mr^{3}}\pa_{3} (\theta e_{3})\cdot\pa_{3} u \md x\nonumber
\\&:= I_{1}+I_{2}
\end{align}
 Thanks to the inequality (\ref{2.1})  and Gagliardo-Nirenberg inequality, $I_{1}$  can be estimated by
 \begin{align}\label{2.5}
I_{1}=&\int_{\mr^{3}}-\pa_{3}u\cdot\nabla u\cdot\pa_{3} u\md x\leq\int_{\mr^{3}}|\pa_{3} u||\nabla u||\pa_{3} u|\md x\nonumber
\\ \leq&
C\|\pa_{3}u\|_{L^{4}}^{2}\|\nabla u\|_{L^{2}}\nonumber
\\
\leq& C\| \pa_{3}u\|_{\dot{B}_{\infty,\infty}^{-r}}\|\pa_{3} u\|_{\dot{H}^{r}}\|\nabla u\|_{L^{2}}\nonumber
\\
\leq& C\| \pa_{3}u\|_{\dot{B}_{\infty,\infty}^{-r}}\|\pa_{3} u\|_{L^{2}}^{1-r}\|\pa_{3}\nabla u\|_{L^{2}}^{r}\|\nabla u\|_{L^{2}}\nonumber
\\
\leq& C\| \pa_{3}u\|_{\dot{B}_{\infty,\infty}^{-r}}^{2}\|\pa_{3}  u\|_{L^{2}}^{2(1-r)}\|\pa_{3}\nabla u\|_{L^{2}}^{2r}+\frac{1}{6}\|\nabla u\|_{L^{2}}^{2}\nonumber
\\
\leq& C\| \pa_{3}u\|_{\dot{B}_{\infty,\infty}^{-r}}^{\frac{2}{1-r}}\|\pa_{3} u\|_{L^{2}}^{2}+\frac{1}{6}(\|\pa_{3}\nabla u\|_{L^{2}}^{2}+\|\nabla u\|_{L^{2}}^{2}).
\end{align}
For $I_{2}$, we have
 by  integrating by parts and H$\rm\ddot{o}$lder inequality
 \begin{align}\label{2.6}
I_{3}\leq\int_{\mr^{3}}|\theta||\pa_{3}\pa_{3} u| \md x
\leq
C\|\theta\|_{L^{2}}\|\pa_{3}\nabla u\|_{L^{2}}
 \leq C\|\theta\|_{L^{2}}^{2}+\frac{1}{6}\|\pa_{3}\nabla u\|_{L^{2}}^{2}.
\end{align}
Now combing the estimates (\ref{2.4}), (\ref{2.5}) and (\ref{2.6}),  it follows
 \begin{align}\label{2.7}
\frac{\md}{\md t}\|\pa_{3} u\|_{L^{2}}^{2}+\|\pa_{3}\nabla u\|_{L^{2}}^{2}
\leq C(\| \pa_{3}u\|_{\dot{B}_{\infty,\infty}^{-r}}^{\frac{2}{1-r}}+\|\nabla u\|_{L^{2}}^{2})(1+\|\pa_{3} u\|_{L^{2}}^{2}).
\end{align}
Thus, the above inequality (\ref{2.7}) yields that with Gronwall's  inequality
  \begin{align}\label{2.8}
 &\|\pa_{3}u\|_{L^{2}}^{2}+\int_{0}^{T}\|\pa_{3}\nabla u\|_{L^{2}}^{2}\md t\nonumber\\
 &\leq C(1+\|\pa_{3} u_{0}\|_{L^{2}}^{2})\exp\left(\int_{0}^{T}\| \pa_{3}u\|_{\dot{B}_{\infty,\infty}^{-r}}^{\frac{2}{1-r}}+\|\nabla u\|_{L^{2}}^{2}\md t\right)\leq C.
 \end{align}
The proof of Lemma \ref{lemma 2.3} is completed.

\textbf{Proof of Theorem \ref{thm1}} \ In the following part,  we shall give the $H^1$ estimates of  $u$ and $\theta$. And firstly, we will bound the $H^1$ norm of  $u$ under the Lemma \ref{lemma 2.3}, then get the $H^1$ norm of $\theta$.
Applying $\nabla$ to the first equation of system $(\ref{1.1})$, then  taking inner product with $\nabla u$, it yields
\begin{align}\label{2.9}
\frac{1}{2}\frac{\md}{\md t}\|\nabla u\|_{L^{2}}^{2}+\|\Delta u\|_{L^{2}}^{2}
&=\int_{\mr^{3}}\nabla(u\cdot\nabla u)\cdot\nabla u\md x+\int_{\mr^{3}}\nabla(\theta e_{3})\cdot\nabla u \md x\nonumber
\\&:= J_{1}+J_{2}
\end{align}
 For $J_{1}$, using the H$\rm\ddot{o}$lder's inequality and Gagliardo-Nirenberg's inequality, and employing the Lemma \ref{lemma 2.2} with $p=6, p_{1}=p_{2}=p_{3}=3$, one has
\begin{align}\label{2.10}
J_{1}&=\int_{\mr^{3}}\nabla u\cdot\nabla u\cdot\nabla u\md x\nonumber
\\&\leq
C\| \nabla u\|_{L^{3}}^{3}
 \leq
C\|\nabla u\|_{L^{2}}^{\frac{3}{2}}\|\nabla u\|_{L^{6}}^{\frac{3}{2}}\nonumber
\\ &\leq
C\|\nabla u\|_{L^{2}}^{\frac{3}{2}}\|\pa_{1}\nabla u\|_{L^{2}}^{\frac{1}{2}}\|\pa_{2}\nabla u\|_{L^{2}}^{\frac{1}{2}}\|\pa_{3}\nabla u\|_{L^{2}}^{\frac{1}{2}}\nonumber
\\ &\leq
C\|\nabla u\|_{L^{2}}^{\frac{3}{2}}\|\Delta u\|_{L^{2}}\|\pa_{3}\nabla u\|_{L^{2}}^{\frac{1}{2}}\nonumber
\\ &\leq
C\|\nabla u\|_{L^{2}}\|\pa_{3}\nabla u\|_{L^{2}}\|\nabla u\|_{L^{2}}^{2}+\frac{1}{4}\|\Delta u\|_{L^{2}}^{2}\nonumber
\\ &\leq
C(\|\nabla u\|_{L^{2}}^{2}+\|\pa_{3}\nabla u\|_{L^{2}}^{2})\|\nabla u\|_{L^{2}}^{2}+\frac{1}{4}\|\Delta u\|_{L^{2}}^{2}.
\end{align}
For $J_{2}$, we have
 \begin{align}\label{2.11}
J_{2}\leq\int_{\mr^{3}}|\theta||\nabla\nabla u| \md x
\leq
C\|\theta\|_{L^{2}}\|\nabla\nabla u\|_{L^{2}}
 \leq C\|\theta\|_{L^{2}}^{2}+\frac{1}{4}\|\Delta u\|_{L^{2}}^{2}.
\end{align}
Inserting the estimates (\ref{2.10}) and (\ref{2.11}) into (\ref{2.9}), we obtain
\begin{align}\label{2.12}
\frac{\md}{\md t}\|\nabla u\|_{L^{2}}^{2}+\|\Delta u\|_{L^{2}}^{2}
\leq C(\|\nabla u\|_{L^{2}}^{2}+\|\pa_{3}\nabla u\|_{L^{2}}^{2})\|\nabla u\|_{L^{2}}^{2}+C\|\theta\|_{L^{2}}^{2}.
\end{align}
Therefore, it can be deduced from the Gronwall inequality and Lemma \ref{lemma 2.2}
\begin{align}\label{2.13}
\|\nabla u\|_{L^{2}}^{2}+\int_{0}^{T}\|\Delta u\|_{L^{2}}^{2}\md t
\leq C.
\end{align}

Next, applying $\nabla$ to the second equation of $(\ref{1.1})$, then  taking inner product with $\nabla \theta$, we have
\begin{align}\label{2.14}
\frac{1}{2}\frac{\md}{\md t}\|\nabla\theta\|_{L^{2}}^{2}+\|\Delta \theta\|_{L^{2}}^{2}
=
\int_{\mr^{3}}\nabla(u\cdot\nabla \theta)\cdot\nabla \theta\md x
\end{align}
 With the incompressibility condition,  the H$\rm\ddot{o}$lder's inequality and Gagliardo-Nirenberg's inequality, it follows that by using the Lemma \ref{lemma 2.3}
\begin{align}\label{2.15}
&\int_{\mr^{3}}\nabla(u\cdot\nabla \theta)\cdot\nabla \theta\md x=\int_{\mr^{3}}\nabla u\cdot\nabla \theta\cdot\nabla \theta\md x\nonumber
\\&\leq
C\| \nabla u\|_{L^{2}}\|\nabla \theta\|_{L^{4}}\|\nabla \theta\|_{L^{4}}
 \leq
C\|\nabla u\|_{L^{2}}\|\nabla \theta\|_{L^{2}}^{\frac{1}{2}}\|\Delta \theta\|_{L^{2}}^{\frac{3}{2}}\nonumber
\\ &\leq
C\|\nabla \theta\|_{L^{2}}^{2}+\frac{1}{4}\|\Delta \theta\|_{L^{2}}^{2} .
\end{align}
Combing (\ref{2.14}) and (\ref{2.15}), we have following
\begin{align}\label{2.16}
\frac{\md}{\md t}\|\nabla\theta\|_{L^{2}}^{2}+\|\Delta \theta\|_{L^{2}}^{2}
\leq C\|\nabla \theta\|_{L^{2}}^{2},
\end{align}
which implies the desired result with Gronwall's inequality
\begin{align}\label{2.17}
\|\nabla \theta\|_{L^{2}}^{2}+\int_{0}^{T}\|\Delta \theta\|_{L^{2}}^{2}\md t
\leq C\|\nabla \theta_{0}\|_{L^{2}}^{2}.
\end{align}

The proof of Theorem 1.1 is completed.
%%%%%%%%%%%%%%%%%%%%%%%%%%%%%%%%%%%%%%%%%%%%%%%%%%%%%%%%%%%%%%%%%%%%%%%%%%%%%%%%%%%%%

\end{document}